\documentclass[11pt]{amsart}
\usepackage{amsmath}
\usepackage{latexsym}
\usepackage{amssymb}
\usepackage{amscd}
\input epsf

\usepackage{amsfonts}
\usepackage[all]{xy}
\usepackage{}

\DeclareFontShape{OT1}{cmr}{b}{n}{<12> cmr12}{}

\newtheorem{proposition}{Proposition}[section]

\newtheorem{lemma}[proposition]{Lemma}
\newtheorem{definition}[proposition]{Definition}
\newtheorem{lemma-definition}[proposition]{Lemma-Definition}
\newtheorem{theorem}[proposition]{Theorem}

\newtheorem{corollary}[proposition]{Corollary}

\newtheorem{example}[proposition]{Example}
\newtheorem{remark}[proposition]{Remark}

\newenvironment{dok}{\par\vspace{-5pt}%
\par\noindent\begingroup%
\leftskip=0em\hspace{0em}{\bf Proof.}}%
{\endgroup\hfill$\Box$}

\newcounter{tmp}

\def\db#1{ \bD^b({#1})}

\def\perf#1{{\mathfrak P}{\mathfrak e}{\mathfrak r}{\mathfrak f}({#1})}
\def\dsing#1{ \bD_{\rm Sg}({#1})}

\def\dsinggr#1{ \bD_{\rm Sg}^{\rm gr}({#1})}

\def\dt{DB}
\def\gdt{\operatorname{DGrB}}

\def\gabe{\operatorname{GrPair}}
\def\ove{\overline}
\def\ul{\underline}

\def\Cok{\operatorname{Cok}}

\def\Ho#1,#2,#3,#4{{\operatorname{Hom}}^{#1}_{#2}({#3}\:,\; {#4})}
\def\Ex#1,#2,#3,#4{{\operatorname{Ext}}^{#1}_{#2}({#3}\:,\; {#4})}

\def\lto{\longrightarrow}

\def\D{{\mathcal D}}

\def\N{{\mathcal N}}
\def\O{{\mathcal O}}
\def\R{{\mathcal R}}
\def\L{{\mathcal L}}

\def\T{{\mathcal T}}

\def\cP{{\mathcal P}}

\def\cS{{\mathcal S}}

\def\ZZ{{\mathbb Z}}

\def\bD{{\mathbf D}}
\def\bR{{\mathbf R}}

\def\bL{{\mathbf L}}

\def\ZZ{{\mathbb Z}}
\def\b1{{\mathbf 1}}
\def\kk{{\mathbf k}}

\def\AA{{\mathbb A}}
\def\NN{{\mathbb N}}
\def\ZZ{{\mathbb Z}}

\def\PP{{\mathbb P}}
\def\pp{{\mathbb P}}

\def\Hom{\operatorname{Hom}}
\def\End{\operatorname{End}}
\def\Ext{\operatorname{Ext}}

\def\Coker{\operatorname{Coker}\,}

\def\Spec{\mathbf{Spec}}
\def\Proj{{\mathbf{Proj}}\,}
\def\PROJ{{\mathbb{P}\mathrm{roj}}\,}
\def\bz{{\mathbf 0}}

\def\id{{\operatorname{id}}}

\def\gr{\operatorname{gr}\!}
\def\qgr{\operatorname{qgr}}
\def\tors{\operatorname{tors}\!}
\def\Gr{\operatorname{Gr}\!}
\def\QGr{\operatorname{QGr}}
\def\Tors{\operatorname{Tors}\!}

\def\mod{\operatorname{mod}\!}
\def\Mod{\operatorname{Mod}\!}

\def\op{\circ}
\def\pr{\operatorname{pr}}

\def\proj{\operatorname{proj}\!}
\def\grproj{\operatorname{grproj}\!}
\def\coh{\operatorname{coh}}
\def\Qcoh{\operatorname{Qcoh}}

\def\qQ{\operatorname{Q}}

\def\Te{P_{\sigma}}

\def\bdot{\begin{picture}(4,4)\put(2,3){\circle*{1.5}}\end{picture}}
\def\Ddot{\bdot}

\title[]{Derived Categories of Coherent Sheaves and Triangulated
Categories of Singularities}

\author[]{Dmitri Orlov}

\address{Algebra Section, Steklov Mathematical Institute RAN,
Gubkin str. 8,  Moscow 119991, RUSSIA}

\email{orlov@mi.ras.ru}

\thanks{
This work was done with  partial financial support from the Weyl
Fund, from grant RFFI 05-01-01034, from grant CRDF Award No
RUM1-2661-MO-05. }

\date{}
\dedicatory{
Dedicated to Yuri Ivanovich Manin on the occasion of his 70th birthday
}

\begin{document}

\begin{abstract}
In this paper we establish an equivalence between the category of
graded D-branes of type B in Landau-Ginzburg models with homogeneous
superpotential $W$ and the triangulated category of singularities of
the fiber of $W$ over zero.  The main result is a theorem that shows
that the graded triangulated category of singularities of the cone
over a projective variety is connected via a fully faithful functor
to the bounded derived category of coherent sheaves on the base of
the cone. This implies that the category of graded D-branes of type
B in  Landau-Ginzburg models with homogeneous superpotential $W$ is
connected via a fully faithful functor to the derived category of
coherent sheaves on the projective variety defined by the equation
$W=0.$
\end{abstract}
\maketitle
\section*{Introduction}

With any algebraic variety $X$ one can naturally associate two
triangulated categories: the bounded derived category $\db{\coh(X)}$
of coherent sheaves and the triangulated subcategory $\perf{X}
\subset \db{\coh(X)}$ of perfect complexes on $X.$ If the variety
$X$ is smooth, then these two categories coincide.  For singular
varieties this is no longer true.  In \cite{Tr} we introduced a new
invariant of a variety $X,$  the triangulated category $\dsing{X}$
of the singularities of $X,$  as the quotient of $\db{\coh(X)}$ by
the full subcategory of perfect complexes $\perf{X}.$ The category
$\dsing{X}$ captures many properties of the singularities of $X.$

Similarly we can define a triangulated category of singularities
$\dsing{A}$ for any Noetherian algebra $A.$ We set $\dsing{A} =
\db{\mod-A}/\perf{A},$ where $\db{\mod-A}$ is the bounded derived
category of finitely generated right $A$\!-modules and $\perf{A}$ is
its triangulated subcategory consisting of objects that are
quasi-isomorphic to bounded complexes of projectives. We will again
call $\perf{A}$ the subcategory of perfect complexes, but usually we
will write $\db{\proj-A}$ instead of $\perf{A},$ since this category
can also be identified with the derived category of the exact
category $\proj-A$ of finitely generated right projective
$A$\!-modules (see, e.g., \cite{Ke}).

The investigation of triangulated categories of singularities is not
only connected with a study of singularities but is mainly inspired by
the Homological Mirror Symmetry Conjecture \cite{K}.  More precisely,
the objects of these categories are directly related to D-branes of
type B (B-branes) in Landau-Ginzburg models.  Such models arise as a
mirrors to Fano varieties \cite{HV}.  For Fano varieties one has the
derived categories of coherent sheaves (B-branes) and given a
symplectic form one can propose a suitable Fukaya category (A-branes).
Mirror symmetry should interchange these two classes of D-branes.
Thus, to extend the Homological Mirror Symmetry Conjecture to the Fano
case, one should describe D-branes in Landau-Ginzburg models.

To specify a Landau-Ginzburg model in general one needs to choose a
target space $X,$ and a holomorphic function $W$ on $X$ called a
superpotential.  The B-branes in the Landau-Ginzburg model are
defined as $W$\!-twisted $\ZZ_2$\!-periodic complexes of coherent
sheaves on $X.$ These are chains $\{ \cdots \stackrel{d}{\to} P_{0}
\stackrel{d}{\to} P_{1} \stackrel{d}{\to} P_{0} \stackrel{d}{\to}
P_{1} \stackrel{d}{\to} P_{0} \cdots \}$ of coherent sheaves in
which the composition of differentials is no longer zero, but is
equal to multiplication by $W$ (see, e.g., \cite{KL,Tr,Or1}).  In
the paper \cite{Tr} we analyzed the relationship between the
categories of B-branes in Landau-Ginzburg models and triangulated
categories of singularities. Specifically, we showed that for an
affine $X$ the product of the triangulated categories of
singularities of the singular fibers of $W$ is equivalent to the
category of B-branes of $(X,W).$

In this paper we consider the graded case.  Let $A=\bigoplus_i A_i$
be a graded Noetherian algebra over a field $\kk.$ We can define the
triangulated category of singularities $\dsinggr{A}$ of $A$ as the
quotient $\db{\gr-A}/\db{\grproj-A},$ where $\db{\gr-A}$ is the
bounded derived category of finitely generated graded right
$A$\!-modules and $\db{\grproj-A}$ is its triangulated subcategory
consisting of objects that are isomorphic to bounded complexes of
projectives.

The graded version of the triangulated category of singularities is
closely related to the category of B-branes in Landau-Ginzburg
models $(X, W)$ equipped with an action of the multiplicative group
$\kk^*$ for which $W$ is semi-invariant.  The notion of grading on
D-branes of type B was defined in the papers \cite{HW, Wa}.  In the
presence of a $\kk^*$-action one can construct a category of graded
B-branes in the Landau-Ginzburg model $(X,W)$ (Definition
\ref{grdbr} and Section \ref{grdbrB}).  Now our Theorem \ref{main3}
gives an equivalence between the category of graded B-branes and the
triangulated category of singularities $\dsinggr{A},$ where $A$ is
such that $\Spec(A)$ is the fiber of $W$ over $0.$

This equivalence allows us to compare the category of graded
B-branes and the derived category of coherent sheaves on the
projective variety that is defined by the superpotential $W.$  For
example, suppose $X$ is the affine space $\AA^N$ and $W$ is a
homogeneous polynomial of degree $d.$ Denote by $Y\subset \PP^{N-1}$
the projective hypersurface of degree $d$ that is given by the
equation $W=0.$  If $d=N,$ then the triangulated category of graded
B-branes $\gdt(W)$ is equivalent to the derived category of coherent
sheaves on the Calabi-Yau variety $Y.$ Furthermore, if $d<N$ (i.e.,
$Y$ is a Fano variety), we construct a fully faithful functor from
$\gdt(W)$ to $\db{\coh(Y)},$ and if $d>N$ (i.e., $Y$ is a variety of
general type), we construct a fully faithful functor from
$\db{\coh(Y)}$ to $\gdt(W)$ (see Theorem \ref{main4}).

This result follows from a more general statement for graded
Gorenstein algebras (Theorem \ref{mai}). It gives a relation between
the triangulated category of singularities $\dsinggr{A}$ and the
bounded derived category $\db{\qgr A},$ where $\qgr A$ is the
quotient of the abelian category of graded finitely generated
$A$\!-modules by the subcategory of torsion modules.  More
precisely, for Gorenstein algebras we obtain a fully faithful
functor between $\dsinggr{A}$ and $\db{\qgr A},$ and the direction
of this functor depends on the Gorenstein parameter $a$ of $A.$ In
particular, when the Gorenstein parameter $a$ is equal to zero, we
obtain an equivalence between these categories.  Finally, the famous
theorem of Serre that identifies $\db{\qgr A}$ with
$\db{\coh(\Proj(A))}$ when $A$ is generated by its first component
allows us to apply this result to geometry.

I am grateful to Alexei Bondal, Anton Kapustin, Ludmil Katzarkov,
Alexander Kuznetsov, Tony Pantev, and Johannes Walcher for very
useful discussions.


\section{Triangulated Categories of Singularities for Graded Algebras}

\subsection{Localization in Triangulated Categories and Semiorthogonal
  Decomposition.}

Recall that a triangulated category $\D$ is an additive category
equipped with the following additional data:
\begin{list}{(\alph{tmp})}%
{\usecounter{tmp}}
\item an additive autoequivalence
$[1]: \D\lto\D,$ which is called a translation functor,
\item a class of exact (or distinguished) triangles
$$
X\stackrel{u}{\lto}Y\stackrel{v}{\lto}Z\stackrel{w}{\lto}X[1],
$$
\end{list}
which must satisfy a certain set of axioms (see \cite{Ve}, also
\cite{GM, Ke, Nee}).

A functor $F : {\D} \lto{\D}'$ between two triangulated categories
is called {\sf exact} if it commutes with the translation functors,
i.e., $F\circ[1]\cong[1]\circ F,$ and transforms exact triangles
into exact triangles.

With any pair $\N\subset \D,$ where $\N$ is a full triangulated
subcategory in a triangulated category $\D,$ we can associate the
quotient category $\D/\N.$ To construct it let us denote by
$\Sigma(\N)$ a class of morphisms $s$ in $\D$ fitting into an exact
triangle
$$
X\stackrel{s}{\lto} Y\lto N\lto X[1]
$$ with $N\in \N.$  It can be checked that $\Sigma(\N)$ is a
multiplicative system.  Define the quotient $ \D/\N $ as the
localization $\D[\Sigma(\N)^{-1}]$ (see \cite{GZ, GM, Ve}). It is a
triangulated category. The translation functor on $\D/\N$ is induced
from the translation functor in the category $\D,$ and the exact
triangles in $\D/\N$ are the triangles isomorphic to the images of
exact triangles in $\D.$ The quotient functor $Q:\D\lto \D/\N$
annihilates $\N.$ Moreover, any exact functor $F: \D\lto \D'$
between triangulated categories for which $F(X)\simeq 0$ when $X\in
\N$ factors uniquely through $Q.$ The following lemma is obvious.
\begin{lemma}\label{adjqu}
Let $\N$ and $\N'$ be full triangulated subcategories of
triangulated categories $\D$ and $\D'$ respectively. Let $F: \D\to
\D'$ and $G: \D'\to \D$ be an adjoint pair of exact functors such
that $F(\N)\subset \N'$ and $G(\N')\subset \N.$ Then they induce
functors
$$
\ove{F}:\D/\N\lto\D'/\N', \qquad \ove{G}:\D'/\N'\lto \D/\N,
$$
which are adjoint as well. Moreover, if the functor $F: \D\to \D'$ is
fully faithful, then the functor $\ove{F}:\D/\N\lto\D'/\N'$ is also fully
faithful.
\end{lemma}

Now recall some definitions and facts concerning admissible
subcategories and semiorthogonal decompositions (see \cite{BK, BO}).
Let ${\N\subset\D}$ be a full triangulated subcategory.  The {\sf
right orthogonal} to ${\N}$ is the full subcategory
${\N}^{\perp}\subset {\D}$ consisting of all objects $M$ such that
${\Hom(N, M)}=0$ for any $N\in{\N}.$ The {\sf left orthogonal}
${}^{\perp}{\N}$ is defined analogously.  The orthogonals are also
triangulated subcategories.

\begin{definition}\label{adm}
Let $I:\N\hookrightarrow\D$ be an embedding of a full triangulated
subcategory $\N$ in a triangulated category $\D.$ We say that ${\N}$
is {\sf right admissible} (respectively {\sf left admissible}) if
there is a right (respectively left) adjoint functor $Q:\D\to \N.$ The
subcategory $\N$ will be called admissible if it is right and left
admissible.
\end{definition}
\begin{remark}\label{semad}
{\rm For the subcategory $\N$ the property of being right admissible
is equivalent to requiring that for each $X\in{\D}$ there be an
exact triangle $N\to X\to M,$ with $N\in{\N}, M\in{\N}^{\perp}.$ }
\end{remark}

\begin{lemma}\label{inv}
Let $\N$ be a full triangulated subcategory in a triangulated category
$\D.$ If $\N$ is right (respectively left) admissible, then the
quotient category $\D/\N$ is equivalent to $\N^{\perp}$ (respectively
${}^{\perp}\N$).  Conversely, if the quotient functor $Q:\D\lto\D/\N$
has a left (respectively right) adjoint, then $\D/\N$ is equivalent to
$\N^{\perp}$ (respectively ${}^{\perp}\N$).
\end{lemma}

If $\N\subset\D$ is a right admissible subcategory, then we say that
the category $\D$ has a weak semiorthogonal decomposition
$\left\langle{\N}^{\perp},{\N}\right\rangle.$ Similarly, if
$\N\subset\D$ is a left admissible subcategory, we say that $\D$ has
a weak semiorthogonal decomposition
$\left\langle{\N},{}^{\perp}{\N}\right\rangle.$

\begin{definition}\label{sd}
A sequence of full triangulated subcategories $({\N}_1, \dots,
{\N}_n)$ in a triangulated category ${\D}$ will be called a {\sf
weak semiorthogonal decomposition} of $\D$ if there is a sequence of
left admissible subcategories $\D_1=\N_1\subset
\D_2\subset\cdots\subset \D_n=\D$ such that ${\N}_p$ is left
orthogonal to $\D_{p-1}$ in $\D_p.$  We will write $
{\D}=\left\langle{\N}_1, \dots, {\N}_n\right\rangle.$ If all $N_p$
are admissible in $\D$ then the decomposition $
{\D}=\left\langle{\N}_1, \dots, {\N}_n\right\rangle$ is called {\sf
semiorthogonal}.
\end{definition}

The existence of a semiorthogonal decomposition on a triangulated
category $\D$ clarifies the structure of $\D.$ In the best scenario,
one can hope that $\D$ has a semiorthogonal decomposition
${\D}=\left\langle{\N}_1, \dots, {\N}_n\right\rangle$ in which each
elementary constituent $\N_p$ is as simple as possible, i.e., is
equivalent to the bounded derived category of finite-dimensional
vector spaces.


\begin{definition}\label{exc}
An object $E$ of a $\kk$\!-linear triangulated category ${\T}$ is
called {\sf exceptional} if  ${\Hom}(E, E[p])=0$ when $p\ne 0,$ and
${\Hom}(E, E)=\kk.$ An {\sf  exceptional collection} in ${\T}$ is a
sequence of exceptional objects $(E_0,\dots, E_n)$ satisfying the
semiorthogonality condition ${\Hom}(E_i, E_j[p])=0$ for all $p$ when
$i>j.$
\end{definition}
If a triangulated category $\D$ has an exceptional collection
$(E_0,\dots, E_n)$ that generates the whole of $\D$ then we say that
the collection is {\sf full}.  In this case $\D$ has a
semiorthogonal decomposition with $\N_p=\langle E_p\rangle.$ Since
$E_{p}$ is exceptional, each of these categories is equivalent to the
bounded derived category of finite dimensional vector spaces.  In
this case we write $ \D=\langle E_0,\dots, E_n \rangle.  $

\begin{definition}\label{strong}
An exceptional collection $(E_0,\dots, E_n)$ is called strong if, in addition,
${\Hom}(E_i, E_j[p])=0$ for all  $i$ and $j$ when $p\ne 0.$
\end{definition}

\subsection{Triangulated Categories of Singularities for Algebras.}

Let $A=\bigoplus_{i\ge 0} A_{i}$ be a Noetherian graded
algebra over a field $\kk.$ Denote by $\mod-A$ and $\gr-A$ the
category of finitely generated right modules and the category of
finitely generated graded right modules respectively.  Note that
morphisms in $\gr-A$ are homomorphisms of degree zero.  These
categories are abelian.  We will also use the notation $\Mod-A$ and
$\Gr-A$ for the abelian categories of all right modules and all
graded right modules and we will often omit the prefix ``right''.
Left $A$\!-modules are will be viewed as right $A^{\op}$\!-modules
and $A-B$ bimodules as right $A^{\op}-B$\!-modules, where $A^{\op}$
is the opposite algebra.

The twist functor $(p)$ on the category $\gr-A$ is defined as
follows: it takes a graded module $M=\mathop{\oplus}_{i} M_{i}$ to
the module $M(p)$ for which $M(p)_{i}= M_{p+i}$ and takes a morphism
$f: M\lto N$ to the same morphism viewed as a morphism between the
twisted modules $f(p): M(p)\lto N(p).$

Consider the bounded derived categories $\db{\gr-A}$ and
$\db{\mod-A}.$ They can be endowed with natural structures of
triangulated categories. The categories $\db{\gr-A}$ and $\db{\mod-A}$
have full triangulated subcategories consisting of objects that are
isomorphic to bounded complexes of projectives. These subcategories
can also be considered as the derived categories of the exact
categories of projective modules $\db{\grproj-A}$ and $\db{\proj-A}$
respectively (see, e.g., \cite{Ke}).  They will be called {\sf the
subcategories of perfect complexes}.  Observe also that the category
$\db{\gr-A}$ (respectively $\db{\mod-A}$) is equivalent to the
category $\bD^b_{\gr-A}(\Gr-A)$ (respectively
$\bD^b_{\mod-A}(\Mod-A)$) of complexes of arbitrary modules with
finitely generated cohomologies (see \cite{Il}). We will tacitly use
this equivalence throughout our considerations.

\begin{definition}
We define triangulated categories of singularities $\dsinggr{A}$
and $\dsing{A}$ as the quotient $\db{\gr-A}/\db{\grproj-A}$
and $\db{\mod-A}/\db{\proj-A}$ respectively.
\end{definition}

\begin{remark}{\rm
As in the commutative case \cite{Tr, Or1}, the triangulated categories
of singularities $\dsinggr{A}$ and $\dsing{A}$ will be trivial if $A$
has finite homological dimension. Indeed, in this case any
$A$\!-module has a finite projective resolution, i.e., the
subcategories of perfect complexes coincide with the full bounded
derived categories of finitely generated modules.}
\end{remark}

Homomorphisms of (graded) algebras $f: A\to B$ induce functors
between the associated derived categories of singularities.
Furthermore, if $B$ has a finite Tor-dimension as an $A$\!-module, then
we get the functor $\stackrel{\bL}{\otimes}_A B$ between the
bounded derived categories of finitely generated modules that maps
perfect complexes to
perfect complexes. Therefore, we get functors between
triangulated categories of singularities
$$
\stackrel{\bL}{\otimes}_A B :\dsinggr{A} \lto \dsinggr{B}
\quad
\text{and}
\quad
\stackrel{\bL}{\otimes}_A B :\dsing{A} \lto \dsing{B}.
$$ If, in addition, $B$ is finitely generated as an $A$\!-module, then
these functors have right adjoints induced from the functor
that sends a complex of $B$\!-modules to itself considered as a
complex of $A$\!-modules.

More generally, suppose $ {}^{}_A{\ul{M}}^{\Ddot}_B$ is a complex of
graded $A-B$ bimodules that as a complex of graded $B$\!-modules is
quasi-isomorphic to a perfect complex.  Suppose that ${}^{}_A
\ul{M}^{\Ddot}$ has a finite Tor-amplitude as a left $A$\!-module.
Then we can define the derived tensor product functor $
\stackrel{\bL}{\otimes}^{}_A\!\ul{M}^{\Ddot}_B :\db{\gr-A}\lto
\db{\gr-B}.$ Moreover, since $\ul{M}^{\Ddot}_B$ is perfect over $B,$
this functor sends perfect complexes to perfect complexes.
Therefore, we get an exact functor
$$
\stackrel{\bL}{\otimes}^{}_A\! \ul{M}^{\Ddot}_B :\dsinggr{A} \lto \dsinggr{B}.
$$
%
%
In the ungraded case
we also get the functor
$
\stackrel{\bL}{\otimes}^{}_A\!\ul{M}^{\Ddot}_B :\dsing{A} \lto \dsing{B}.
$

\subsection{Morphisms in Categories of Singularities.}

In general, it is not easy to calculate spaces of morphisms between
objects in a quotient category.  The following lemma and proposition
provide some information about the morphism spaces in triangulated
categories of singularities.

\begin{lemma}\label{lfint}
For any object $T\in\dsinggr{A}$ (respectively $T\in\dsing{A}$) and
for any sufficiently large $k,$ there is a module $M\in \gr-A$
(respectively $M\in\mod-A$) depending on $T$ and $k$ and such that
$T$ is isomorphic to the image of $M[k]$ in the triangulated
category of singularities.  If, in addition, the algebra $A$ has
finite injective dimension, then for any sufficiently large $k$ the
corresponding module $M$ satisfies $\Ext^{i}_A(M, A)=0$ for all
$i>0.$
\end{lemma}
\begin{dok} The object $T$ is represented by a bounded complex of modules
$\ul{T}^{\Ddot}.$ Choose a bounded above projective resolution
$\ul{P}^{\Ddot}\stackrel{\sim}{\to} \ul{T}^{\Ddot}$ and a
  sufficiently large $k\gg 0.$ Consider
the stupid truncation $\sigma^{\ge -k+1}\ul{P}^{\Ddot}$ of
$\ul{P}^{\Ddot}.$  Denote by $M$ the cohomology module
$H^{-k+1}(\sigma^{\ge -k+1}\ul{P}^{\Ddot}).$  Clearly $T\cong M[k]$
in $\dsinggr{A}$ (respectively $\dsing{A}).$

If now $A$ has finite injective dimension, then morphism spaces
$\Hom(\ul{T}^{\cdot}, A[i])$ in $\db{\gr-A}$ (respectively
$\db{\mod-A}$) are trivial for all but finitely many $i\in \ZZ.$ So
if  $M$ corresponds to $T$ and a sufficiently large $k,$ then we
will have $\Ext^i_A(M, A)=0$ for all $i>0.$
\end{dok}

\

\begin{proposition}\label{stab}
Let $M$  be an $A$\!-module such that $\Ext^{i}_{A}(M, A)=0$ for all
$i>0.$ Then for any $A$\!-module $N$ we have
$$
\Hom_{\dsing{A}}(M, N)\cong \Hom_{A}(M, N)/\R,
$$
where $\R$ is the subspace of elements factoring through a projective module,
i.e.,  $e\in \R$ iff $e=\beta\alpha$ with $\alpha: M\to P$
and $\beta: P\to N,$ where $P$ is projective.
If $M$ is a graded module, then for any graded $A$\!-module $N$
$$
\Hom_{\dsinggr{A}}(M, N)\cong \Hom_{\gr-A}(M, N)/\R.
$$
\end{proposition}
\begin{dok}
We will  discuss only the graded case.  By the definition of localization any
morphism from $M$ to $N$ in $\dsinggr{A}$ can be represented by a pair
\begin{equation}\label{domik}
M\stackrel{a}{\lto}\ul{T}^{\Ddot}\stackrel{s}{\longleftarrow} N
\end{equation}
of morphisms in $\db{\gr-A}$ such that the cone $\ul{C}^{\Ddot}(s)$
is a perfect complex.  Consider a bounded above projective
resolution $\ul{Q}^{\Ddot}\to N$ and its stupid truncation
$\sigma^{\ge -k} \ul{Q}^{\Ddot}$ for sufficiently large $k.$ There
is an exact triangle
$$
E[k]\lto \sigma^{\ge -k} \ul{Q}^{\Ddot}\lto N\stackrel{s'}{\lto} E[k+1],
$$ where $E$ denotes the module $H^{-k}(\sigma^{\ge
-k}\ul{Q}^{\Ddot}).$  Choosing $k$ to be sufficiently large, we can
guarantee that $\Hom(\ul{C}^{\Ddot}(s), E[i])=0$ for all $i>k.$
Using the triangle
$$
\ul{C}^{\Ddot}(s)[-1]\lto N \stackrel{s}{\lto} \ul{T}^{\Ddot}\lto
\ul{C}^{\Ddot}(s),
$$
we find that the map $s': N\to E[k+1]$ can be lifted to a map
$\ul{T}^{\Ddot}\to E[k+1].$ The  map $\ul{T}^{\Ddot}\to E[k+1]$
induces a pair of the form
\begin{equation}\label{domik2}
M\stackrel{a'}{\lto}E[k+1]\stackrel{s'}{\longleftarrow} N,
\end{equation}
and this pair gives the same morphism in $\dsing{A}$ as the pair
(\ref{domik}). Since $\Ext^{i}(M, P)=0$ for all $i>0$ and any
projective module $P,$ we obtain
$$
\Hom(M,\; (\sigma^{\ge -k}\ul{Q}^{\Ddot})[1])=0.
$$
Hence, the map $a':M\to E[k+1]$ can be lifted to a map $f$  that
completes the  diagram
$$
\xymatrix{
M\ar[rr]^-{f}\ar[dr]_-{a'}
&& N\ar[dl]^-{s'}
\\
&E[k+1]}
$$
Thus, the map $f$ is equivalent to the map (\ref{domik2}) and, as a
consequence, to the map (\ref{domik}). Hence, any morphism from $M$
to $N$ in $\dsinggr{A}$ is represented by  a morphism from $M$ to
$N$ in the category $\db{\gr-A}.$

Now if $f$ is the $0$\!-morphism in $\dsinggr{A},$ then without
loss of generality we can assume that the map $a$ is the zero map.
In this case we will have $a'=0$ as well.  This implies that $f$
factors through a morphism $M\to \sigma^{\ge -k}\ul{Q}^{\Ddot}.$ By
the assumption on $M,$ any such morphism can be lifted to a morphism
$M\to Q^0.$ Hence, if $f$ is the $0$\!-morphism in $\dsinggr{A},$
then it factors through $Q^0.$ The same proof works in the ungraded case
(see \cite{Tr}).
\end{dok}

\

\bigskip

\noindent Next we describe a useful construction utilizing the
previous statements. Let $\ul{M}^{\Ddot}$ and $\ul{N}^{\Ddot}$ be
two bounded complexes of (graded) $A$\!-modules. Assume that $\Hom
(\ul{M}^{\Ddot}, A[i])$ in the bounded derived categories of
$A$\!-modules are trivial except for finite number of $i\in \ZZ.$ By
Lemma \ref{lfint}, for sufficiently large $k$ there are modules $M,
N\in \gr-A$ (resp. $M, N\in\mod-A$) such that $\ul{M}^{\Ddot}$ and
$\ul{N}^{\Ddot}$ are isomorphic to the images of $M[k]$ and $N[k]$
in the triangulated category of singularities.  Moreover, it follows
immediately from the assumption on $\ul{M}^{\Ddot}$ and the
construction of $M$ that for any sufficiently large $k$ we have
$\Ext^i_A(M, A)=0$ whenever $i>0.$ Hence, by Proposition \ref{stab},
we get
$$
\Hom_{\dsinggr{A}}(\ul{M}^{\Ddot}, \ul{N}^{\Ddot})\cong
\Hom_{\dsinggr{A}}(M, N)\cong
\Hom_{A}(M, N)/\R,
$$
where $\R$
is the subspace of elements factoring through a projective module.
This procedure works in the ungraded situation as well.

\section{Categories of Coherent Sheaves and Categories of
  Singularities}

\subsection{Quotient Categories of Graded Modules.}

Let $A=\bigoplus_{i\ge 0} A_{i}$ be a Noetherian graded
algebra.  We suppose that $A$ is connected, i.e., $A_0=\kk.$ Denote
by $\tors-A$ the full subcategory of $\gr-A,$ which consists of all
graded $A$\!--modules that are finite dimensional over $\kk.$

An important role will be played by the quotient abelian category
$\qgr A=\gr-A/\tors-A.$ It has the following explicit description.
The objects of $\qgr A$ are the objects of $\gr-A$ (we denote by
$\pi{M}$ the object in $\qgr A$ that corresponds to a module
$M$). The morphisms in $\qgr A$ are given by
\begin{equation}
\Hom_{\qgr}(\pi{M}, \pi{N}):=\lim_{\stackrel{\lto}{M'}}
\Hom_{\gr\;}(M', N),
\end{equation}
where $M'$ runs over submodules of $M$ such that $M/M'$ is finite
dimensional.

Given a graded $A$\!-module $M$ and an integer $p,$ the graded
$A$\!-submodule $\bigoplus_{i\ge p} M_i$ of $M$ is denoted by
$M_{\ge p}$ and is called the $p$-th tail of $M.$ In the same way, we
can define the $p$\!-th tail $\ul{M}^{\Ddot}_{\ge p}$  of any complex
of modules $\ul{M}^{\Ddot}.$ Since $A$ is Noetherian, we have
$$
\Hom_{\qgr}(\pi{M}, \pi{N})=\lim_{p\to \infty} \Hom_{\gr\;}(M_{\ge p},
N).
$$
We will also identify $M_p$ with the quotient $M_{\ge p}/M_{\ge p+1}.$

Similarly, we can consider the subcategory $\Tors-A\subset\Gr-A$ of
torsion modules.  Recall that a module $M$ is called torsion if for
any element $x\in M$ one has $x A_{\ge p}=0$ for some $p.$ Denote by
$\QGr A$ the quotient category $\Gr-A/\Tors-A.$ The category $\QGr A$
contains $\qgr A$ as a full subcategory.  Sometimes it is convenient
to work in $\QGr A$ instead of $\qgr A.$

Denote by $\varPi$ and $\pi$ the canonical projections of $\Gr-A$ to
$\QGr A$ and of $\gr-A$ to $\qgr A$ respectively. The functor $\varPi$
has a right adjoint $\varOmega,$ and moreover, for any $N\in\Gr-A,$
\begin{equation}\label{funad}
\varOmega\varPi N\cong \bigoplus_{n=-\infty}^{\infty}
\Hom_{\QGr}(\varPi A, \varPi N(n)).
\end{equation}
For any $i\in\ZZ$ we can consider the
full abelian subcategories $\Gr-A_{\ge i} \subset \Gr A$ and
$\gr-A_{\ge i} \subset \gr A,$
which consist of all modules $M$ such that $M_{p}=0$ when $p<i.$
The natural projection functor $\varPi_i : \Gr-A_{\ge i}\lto \QGr-A$
has a right adjoint $\varOmega_i$ satisfying
$$
\varOmega_i\varPi_i N\cong \bigoplus_{n=i}^{\infty} \Hom_{\QGr}(\varPi
A, \varPi_i N(n)).
$$
Since the category $\QGr A$ is an abelian category with enough injectives,
there is a right derived functor
$$
\bR \varOmega_i: \bD^{+}(\QGr A)\lto \bD^{+}(\Gr-A_{\ge i})
$$
defined as
\begin{equation}\label{adjback}
\bR \varOmega_i M\cong \bigoplus_{k=i}^{\infty} \bR\Hom_{\QGr}(\varPi
A, M(k)).
\end{equation}

Assume now that the algebra $A$ satisfies condition ``$\chi$'' from
\cite[Sec.~3]{AZ}.  We recall that by definition, a connected
Noetherian graded algebra $A$ satisfies condition ``$\chi$'' if for
every $M\in\gr-A$ the grading on the space $\Ext_A^i(\kk, M)$ is
right bounded for all $i.$  In this case it was proved in
\cite[Prop.~3.14]{AZ}
 that the restrictions of the functors $\varOmega_i$
to the subcategory $\qgr A$ give functors $\omega_i:\qgr A\lto
\gr-A_{\ge i}$ that are right adjoint to $\pi_i.$ Moreover, it follows
from \cite[Th.~7.4]{AZ} that the functor $\omega_i$ has a right
derived
$$
\bR \omega_i: \bD^{+}(\qgr A)\lto \bD^{+}(\gr-A_{\ge i})
$$
and all $\bR^{j}\omega_i\in \tors-A$ for $j>0.$

If, in addition, the algebra $A$ is Gorenstein (i.e.,  if it has a
finite injective dimension $n$ and $D(\kk)=\bR \Hom_A(\kk, A)$ is
isomorphic to $\kk(a)[-n]$), we obtain the right derived functor
$$
\bR \omega_i: \db{\qgr A}\lto \db{\gr-A_{\ge i}}
$$
between bounded derived categories (see \cite[Cor.~4.3]{YZ}). It is
important to note that the functor $\bR\omega_i$ is fully faithful
because $\pi_i\bR\omega_i$ is isomorphic to the identity functor
(\cite[Prop.~7.2]{AZ}).

\subsection{Triangulated Categories of Singularities for Gorenstein
  Algebras.}

The main goal of this section is to establish a connection between
the triangulated category of singularities $\dsinggr{A}$ and the
derived category $\db{\qgr A},$ in the case of a Gorenstein algebra
$A.$

When the algebra $A$ has finite injective dimension
as both a right and  left module over itself
(i.e., $A$ is a dualizing complex for itself)
 we get two functors
\begin{align}
&D^{\hphantom{\op}}:=\bR \Hom_{A^{\hphantom{\op}}}(-, A):
  \db{\gr-A^{\hphantom{\op}}}^{\op}\lto
  \db{\gr-A^{\op}},\label{dualD}\\
&D^{\op}:= \bR\Hom_{A^{\op}}(-, A):
\db{\gr-A^{\op}}^{\op}\lto \db{\gr-A^{\hphantom{\op}}},
\end{align}
which are quasi-inverse triangulated equivalences (see
\cite[Prop.~3.5]{Ye}).

\begin{definition}
We say that a connected graded Noetherian algebra $A$ is Gorenstein
if it has a finite injective dimension $n$ and $D(\kk)=\bR
\Hom_A(\kk, A)$ is isomorphic to $\kk(a)[-n]$ for some integer $a,$
which is called the Gorenstein parameter of $A.$ (Such an algebra is
also called AS-Gorenstein, where ``AS'' stands for ``Artin-Schelter''.)
\end{definition}

\begin{remark}\label{gorhi}{\rm
It is known (see \cite[Cor.~4.3]{YZ}) that any Gorenstein algebra
satisfies condition ``$\chi$'' and  for any Gorenstein algebra $A$ and
for any $i\in \ZZ$ we have derived functors
$$
\bR \omega_i: \db{\qgr A}\lto \db{\gr-A_{\ge i}}
$$
that are fully faithful.
}
\end{remark}

Now we describe the images of the functors $\bR \omega_i.$ Denote by
$\D_i$ the subcategories of $\db{\gr-A}$ that are the images of the
composition of $\bR\omega_i$ and the natural inclusion of
$\db{\gr-A_{\ge i}}$ to $\db{\gr-A}.$ All $\D_i$ are equivalent to
$\db{\qgr-A}.$ Further, for any integer $i$ denote by $\cS_{< i}(A)$
(or simply $\cS_{<i}$) the full triangulated subcategory of
$\db{\gr-A}$ generated by the modules $\kk(e)$ with $e>-i.$ In other
words, the objects of $\cS_{<i}$ are complexes $\ul{M}^{\Ddot}$ for
which the tail $\ul{M}^{\Ddot}_{\ge i}$ is isomorphic to
zero. Analogously, we define $\cS_{\ge i}$ as the triangulated
subcategory that is generated by the objects $\kk(e)$ with $e\le -i.$
In other words, the objects of $\cS_{\ge i}$ are complexes of torsion
modules from $\gr-A_{\ge i}.$ It is clear that $\cS_{<i}\cong \cS_{<
0}(-i)$ and $\cS_{\ge i}\cong \cS_{\ge 0}(-i).$

Furthermore, denote by $\cP_{< i}$ the full triangulated subcategory
of $\db{\gr-A}$ generated by the free modules $A(e)$ with $e>-i$ and
denote by $\cP_{\ge i}$ the triangulated subcategory that is
generated by the free modules $A(e)$ with $e\le -i.$ As above, we have
$\cP_{<i}\cong \cP_{< 0}(-i)$ and $\cP_{\ge i}\cong \cP_{\ge 0}(-i).$

\begin{lemma}\label{admsub}
Let $A=\bigoplus_{i\ge 0} A_{i}$  be a connected graded
Noetherian   algebra. Then  the subcategories  $\cS_{< i}$ and
$\cP_{<i}$ are left and respectively right admissible for any $i\in
\ZZ.$ Moreover, there are  weak semiorthogonal decompositions
\begin{align}
&\db{\gr-A}=\langle \cS_{< i}, \db{\gr-A_{\ge i}}\rangle,
\qquad
\db{\tors-A}=\langle \cS_{< i}, \cS_{\ge i} \rangle,\label{first}\\
&\db{\gr-A}=\langle  \db{\gr-A_{\ge i}}, \cP_{<i} \rangle,
\qquad
\db{\grproj-A}=\langle  \cP_{\ge i}, \cP_{<i} \rangle.\label{second}
\end{align}
\end{lemma}
\begin{dok}
For any complex $\ul{M}^{\Ddot}\in\db{\mod-A}$ there is an exact triangle
of the form
$$
\ul{M}^{\Ddot}_{\ge i}\lto \ul{M}^{\Ddot}\lto
\ul{M}^{\Ddot}/\ul{M}^{\Ddot}_{\ge i}.
$$
By definition, the object $\ul{M}^{\Ddot}/\ul{M}^{\Ddot}_{\ge i}$
belongs to $\cS_{< i},$ and  the object
$\ul{M}^{\Ddot}_{\ge i}$ is in the left orthogonal ${}^{\perp}\cS_{< i}.$
Hence, by Remark \ref{semad}, $\cS_{< i}$ is left admissible.
Moreover, $\ul{M}^{\Ddot}_{\ge i}$ also belongs to $\db{\gr-A_{\ge i}},$
i.e., $\db{\gr-A_{\ge i}}\cong {}^{\perp}\cS_{< i}$ in the category
$\db{\gr-A}.$
If $\ul{M}^{\Ddot}$ is a complex of torsion modules, then
$\ul{M}^{\Ddot}_{\ge i}$ belongs to $\cS_{\ge i}.$
Thus, we obtain both  decompositions of (\ref{first}).

To prove the existence of the decompositions (\ref{second}) we first
note that due to the connectedness of $A,$ any finitely generated
graded projective $A$\!-module is free.  Second, any finitely
generated free module $P$ has a canonical split decomposition of the
form
$$
0\lto P_{<i}\lto P\lto P_{\ge i}\lto 0,
$$ where $P_{<i}\in \cP_{<i}$ and $P_{\ge i} \in \cP_{\ge i}.$ Third,
any bounded complex of finitely generated $A$\!-modules
$\ul{M}^{\Ddot}$ has a bounded above free resolution
$\ul{P}^{\Ddot}\to \ul{M}^{\Ddot}$ such that $P^{-k}\in\cP_{\ge
i}$ for all $k\gg 0.$ This implies that the object
$\ul{P}^{\Ddot}_{<i}\in\cP_{<i}$ from the exact sequence of
complexes
$$
0\lto \ul{P}^{\Ddot}_{<i}\lto \ul{P}^{\Ddot}\lto
\ul{P}^{\Ddot}_{\ge i}\lto 0
$$
is a bounded complex.  Since $\ul{P}^{\Ddot}$ is quasi-isomorphic
to a bounded complex, the complex $\ul{P}^{\Ddot}_{\ge i}$ is also
quasi-isomorphic to some bounded complex $\ul{K}^{\Ddot}$ from
$\db{\gr-A_{\ge i}}.$ Thus, any object $\ul{M}^{\Ddot}$ has a
decomposition
$$
\ul{P}^{\Ddot}_{<i}\lto \ul{M}^{\Ddot}\lto \ul{K}^{\Ddot},
$$ where $\ul{P}^{\Ddot}_{<i}\in\cP_{<i}$ and
$\ul{K}^{\Ddot}\in\db{\gr-A_{\ge i}}.$ This proves the decompositions
(\ref{second}).
\end{dok}

\begin{lemma}\label{dec} Let $A=\bigoplus_{i\ge 0} A_{i}$  be a
  connected graded Noetherian algebra
that is Gorenstein. Then the subcategories  $\cS_{\ge i}$ and $\cP_{\ge i}$
are  right and respectively left admissible.
Moreover, for any $i\in \ZZ$
there are  weak semiorthogonal decompositions
\begin{equation}\label{third}
\db{\gr-A_{\ge i}}=\langle \D_i, \cS_{\ge i}\rangle,
\quad
\db{\gr-A_{\ge i}}=\langle \cP_{\ge i}, \T_i \rangle,
\end{equation}
where the subcategory $\D_i$ is equivalent to $\db{\qgr A}$ under
the functor $\bR\omega_i,$ and $\T_i$ is equivalent to $\dsinggr{A}.$
\end{lemma}
\begin{dok}
The functor $\bR\omega_i$ is fully faithful and has the left adjoint
$\pi_i.$
Thus, we obtain a semiorthogonal decomposition
$$
\db{\gr-A_{\ge i}}=\langle \D_i, {}^{\perp}\D_i\rangle,
$$
where $\D_i\cong\db{\qgr A}.$
Furthermore, the orthogonal ${}^{\perp}\D_i$ consists of  all objects
$\ul{M}^{\Ddot}$
satisfying $\pi_i(\ul{M}^{\Ddot})=0.$ Thus, ${}^{\perp}\D_i$ coincides
with $\cS_{\ge i}.$ Hence, $\cS_{\ge i}$ is right admissible in
$\db{\gr-A_{\ge i}},$ which is right admissible in all of $\db{\gr-A}.$
This implies that $\cS_{\ge i}$ is right admissible in $\db{\gr-A}$ as
well.

The functor $D$ from (\ref{dualD}) establishes an equivalence of the
subcategory $\cP_{\ge i}(A)^{\op}$ with the subcategory
$\cP_{<-i+1}(A^{\op}),$ which is right admissible by Lemma
\ref{admsub}.  Hence, $\cP_{\ge i}(A)$ is left admissible and there is
a decomposition of the form
$$
\db{\gr-A_{\ge i}}=\langle \cP_{\ge i}, \T_i \rangle
$$
with some $\T_i.$

Now applying Lemma \ref{adjqu}, to the full embedding of $\db{\gr-A_{\ge i}}$
to $\db{\gr-A}$ and using Lemma \ref{inv} we get a fully faithful
functor from $\T_i\cong\db{\gr-A_{\ge i}}/\cP_{\ge i}$ to
$\dsinggr{A}=\db{\gr-A}/\db{\grproj-A}.$
Finally, since this functor is essentially surjective on objects,
it is actually an equivalence.
\end{dok}

\begin{theorem}\label{mai}
Let $A=\bigoplus_{i\ge 0} A_{i}$ be a connected graded
Noetherian algebra that is Gorenstein with Gorenstein parameter
$a.$ Then the triangulated categories $\dsinggr{A}$ and $\db{\qgr
A}$ are related as follows:
\begin{enumerate}
\item[{\rm (i)}]
if $a>0,$
there are  fully faithful functors
$\Phi_i: \dsinggr{A}\lto \db{\qgr A}$ and   semiorthogonal decompositions
$$
\db{\qgr A}=\langle \pi A(-i-a+1),\dots, \pi A(-i), \Phi_i \dsinggr{A}\rangle,
$$
where $\pi:\db{\gr-A}\lto\db{\qgr A}$ is the natural projection;
\item[{\rm (ii)}]
if $a<0,$
there are  fully faithful functors
$\Psi_i: \db{\qgr A}\lto \dsinggr{A}$ and   semiorthogonal decompositions
$$
\dsinggr{A}=\langle q\kk(-i),\dots, q\kk (-i+a+1), \Psi_i \db{\qgr A}\rangle,
$$
where $q:\db{\gr-A}\lto\dsinggr{A}$ is the natural projection;
\item[{\rm (iii)}]
if $a=0,$  there is an equivalence
$\dsinggr{A}\stackrel{\sim}{\lto} \db{\qgr A}.$
\end{enumerate}
\end{theorem}
\begin{dok}
Lemmas \ref{admsub} and \ref{dec} gives us that the subcategory $\T_i$
is admissible in $\db{\gr-A}$ and the right orthogonal $\T_i^{\perp}$
has a weak semiorthogonal decomposition of the form
\begin{equation}\label{ti}
\T_i^{\perp}=\langle \cS_{<i}, \cP_{\ge i}\rangle.
\end{equation}
Now let us describe the right orthogonal to the subcategory $\D_i.$
First, since $A$ is Gorenstein the functor $D$ takes the subcategory
$\cS_{\ge i}(A)$ to the subcategory $\cS_{<-i-a+1}(A^{\op}).$ Hence,
$D$ sends the right orthogonal $\cS_{\ge i}^{\perp}(A)$ to the left
orthogonal ${}^{\perp}\cS_{<-i-a+1}(A^{\op})$ which coincides with the
right orthogonal $\cP_{<-i-a+1}^{\perp}(A^{\op})$ by Lemma
\ref{admsub}.  Therefore, the subcategory $\cS_{\ge i}^{\perp}$
coincides with ${}^{\perp}\cP_{\ge i+a}.$ On the other hand, by Lemmas
\ref{admsub} and \ref{dec} we have that
$$
{}^{\perp}\cP_{\ge i+a}=\cS_{\ge i}^{\perp}\cong\langle \cS_{< i},
\D_i \rangle.
$$
This implies that the right orthogonal $\D_{i}^{\perp}$ has the
following decomposition
\begin{equation}\label{di}
\D_i^{\perp}=\langle \cP_{\ge i+a},  \cS_{< i} \rangle.
\end{equation}

Assume that $a\ge 0.$ In this case, the decomposition (\ref{di})
is not only semiorthogonal, but is in fact mutually orthogonal, because
$\cP_{\ge i+a}\subset\db{\gr-A_{\ge i}}.$ Hence, we can interchange
$\cP_{\ge i+a}$ and $\cS_{< i},$ i.e.
$$
\D_i^{\perp}=\langle \cS_{< i}, \cP_{\ge i+a} \rangle.
$$
Thus, we obtain that $\D_i^{\perp}\subset\T_i^{\perp}$ and,
consequently, $\T_i$ is a full subcategory of $\D_i.$ Moreover, we can
describe the right orthogonal to $\T_i$ in $\D_i.$ In fact, there is
a decomposition
$$
\cP_{\ge i}=\langle \cP_{\ge i+a}, \cP_i^{a}\rangle,
$$
where $\cP_i^{a}$ is the subcategory generated by the modules
$A(-i-a+1),\dots, A(-i).$ Moreover, these modules form an exceptional
collection. Thus, the category $\D_i$ has the semiorthogonal
decomposition
$$
\D_i=\langle  A(-i-a+1),\dots,  A(-i), \T_i\rangle.
$$
Since $\D_i\cong \db{\qgr A}$ and $\T_i\cong \dsinggr{A},$
we obtain the decomposition
$$
\db{\qgr A}\cong\langle \pi A(-i-a+1),\dots, \pi A(-i), \Phi_i
\dsinggr{A}\rangle,
$$
where the fully faithful functor $\Phi_i$ is the composition
$\dsinggr{A}\stackrel{\sim}{\to}
\T_i\hookrightarrow\db{\gr-A}\stackrel{\pi}{\to}
\db{\qgr A}.
$

Assume now that $a\le 0.$ In this case, the decomposition (\ref{ti})
is not only semiorthogonal but is in fact mutually orthogonal, because
the algebra $A$ is Gorenstein and $\bR\Hom_A(\kk, A)=\kk(a)[-n]$ with
$a\le 0.$ Hence, we can interchange $\cP_{\ge i}$ and $\cS_{< i},$
i.e.,
$$
\T_i^{\perp}=\langle \cP_{\ge i}, \cS_{< i} \rangle.
$$
Now we see that $\T_i^{\perp}\subset\D_{i-a}^{\perp},$ and
consequently, $\D_{i-a}$ is the full subcategory of $\T_i.$ Moreover,
we can describe the right orthogonal to $\D_{i-a}$ in $\T_i.$
In fact, there is a decomposition of the form
$$
\cS_{< i-a}=\langle \cS_{< i}, \kk(-i),\dots, \kk(-i+a+1)\rangle.
$$
Therefore, the category $\T_i\cong \dsinggr{A}$ has a semiorthogonal
decomposition
of the form
\begin{equation}\label{tidec}
\T_i=\langle  \kk(-i),\dots,  \kk(-i+a+1), \D_{i-a}\rangle.
\end{equation}
Since $\D_{i-a}\cong \db{\qgr A}$ and $\T_i\cong \dsinggr{A},$
we obtain the decomposition
$$
\dsinggr{A}\cong \langle  q\kk(-i),\dots,  q\kk(-i+a+1),
\Psi_{i}\db{\qgr A}\rangle,
$$
where the fully faithful functor $\Psi_i$ can be defined as the
composition
$\db{\qgr A}\stackrel{\sim}{\to}
\D_{i-a}\hookrightarrow\db{\gr-A}\stackrel{q}{\to}
\dsinggr{A}.
$
If $a=0,$ then we get equivalence.
\end{dok}

\begin{remark}{\rm
It follows from the construction that the functor $\Psi_{i+a}$ from the bounded derived category
$\db{\qgr A}$ to $\dsinggr{A}$ is the composition of the functor
$\bR\omega_i: \db{\qgr A} \lto \db{\gr-A_{\ge i}},$ which is given by formula
(\ref{adjback}), the natural embedding $\db{\gr-A_{\ge i}}\hookrightarrow \db{\gr-A},$
and the projection $\db{\gr-A}\stackrel{q}{\to}
\dsinggr{A}.$
}
\end{remark}

Let us consider two limiting cases. The first case is that the
algebra $A$ has finite homological dimension. In this case the
triangulated category of singularities $\dsinggr{A}$ is trivial, and
hence the Gorenstein parameter $a$ is nonnegative and the derived
category $\db{\qgr A}$ has a full exceptional collection
$\sigma=\big( \pi A(0),\dots, \pi A(a-1)\big).$ More precisely, we
have the following:

\begin{corollary}
Let $A=\bigoplus_{i\ge 0} A_{i}$  be a connected graded
Noetherian algebra that is Gorenstein with Gorenstein parameter
$a.$ Suppose that $A$ has finite homological dimension. Then, $a\ge
0$ and the derived category $\db{\qgr A}$ has a full strong
exceptional collection $ \sigma=\left(\pi A(0),\dots, \pi
A(a-1)\right). $ Moreover, the  category $\db{\qgr A}$ is equivalent
to the derived category $\db{\mod-\qQ(A)}$ of finite (right) modules
over the algebra $\qQ(A)
:=\End_{\gr-A}\left(\bigoplus\limits_{i=0}^{a-1} A(i)\right)$ of
homomorphisms of $\sigma.$
\end{corollary}
\begin{dok}
Since $A$ has a finite homological dimension, the category
$\dsinggr{A}$ is trivial. By Theorem \ref{mai} we get that $a\ge
0$ and that $\db{\qgr A}$ has a full exceptional collection $
\sigma=\left( \pi A(0),\dots, \pi A(a-1)\right). $ Consider the
object $\Te=\bigoplus\limits_{i=0}^{a-1} \pi A(i)$ and  the
functor
$$
\Hom(\Te, -):\qgr A\lto \mod-\qQ(A),
$$
where $\qQ(A)=\End_{\qgr A}\left(\bigoplus\limits_{i=0}^{a-1} \pi
A(i)\right)= \End_{\gr-A}\left(\bigoplus\limits_{i=0}^{a-1}
A(i)\right)$ is the algebra of homomorphisms of the exceptional
collection $\sigma.$ It is easy to see that this functor has a
right derived functor
$$
\bR\Hom(\Te, -):\db{\qgr A}\lto \db{\mod-\qQ(A)}
$$
(e.g., as a composition $\bR\omega_0$ and $\Hom \big(
\bigoplus\limits_{i=0}^{a-1} A(i), - \big)$). The standard
reasoning (see, e.g., \cite{Bon} or \cite{BK}) now shows that the
functor $\bR\Hom(\Te, -)$ is an equivalence.
\end{dok}

\begin{example}
{\rm As an application we obtain a well-known result (see \cite{Be})
asserting the existence of a full exceptional collection in the
bounded derived category of coherent sheaves on the projective space
$\PP^n.$ This result follows immediately if we take
$A=\kk[x_0,\dots, x_n]$ with its standard grading. More generally,
if we take $A$ to be the polynomial algebra $\kk[x_0, \ldots,x_n]$
graded by $\deg x_i =a_i,$ then we get a full exceptional collection
$\big(\O,\dots, \O(\sum_{i=0}^{n}a_i-1)\big)$ in the bounded derived
category of coherent sheaves on the weighted projective space
$\PP(a_0,\dots, a_n)$ considered as a smooth orbifold (see \cite{Ba,
AKO}). It is also true for noncommutative (weighted) projective
spaces \cite{AKO}. }
\end{example}

Another limiting case is that the algebra $A$ is finite-dimensional
over the base field (i.e., $A$ is a Frobenius algebra).  In this case
the category $\qgr A$ is trivial, and hence the triangulated
category of singularities $\dsinggr{A}$ has a full exceptional
collection (compare with \cite{Hap} 10.10). More precisely, we get
the following:

\begin{corollary}
Let $A=\bigoplus_{i\ge 0} A_{i}$  be a connected graded
Noetherian algebra that is Gorenstein with Gorenstein parameter
$a.$ Suppose that $A$ is  finite  dimensional over the field $\kk.$
Then $a\le 0,$ and the triangulated category of singularities
$\dsinggr{A}$ has a full exceptional collection $ \left(
q\kk(0),\dots, q\kk (a+1)\right), $ where
$q:\db{\gr-A}\lto\dsinggr{A}$ is the natural projection. Moreover,
the  triangulated category $\dsinggr{A}$ is equivalent to the
derived category $\db{\mod-\qQ(A)}$ of finite (right) modules over
the algebra $\qQ(A)=\End_{\gr-A}\left(\bigoplus\limits_{i=a+1}^{0}
A(i)\right).$
\end{corollary}
\begin{dok}
Since $A$ is finite-dimensional, the derived category $\db{\qgr
A}$ is trivial. By Theorem \ref{mai} we get that $a\le 0$ and
$\dsinggr{A}$ has a full exceptional collection $ \left(
q\kk(0),\dots, q\kk (a+1)\right). $ Unfortunately, this collection
is not strong. However, we can replace it by the dual exceptional
collection which is already strong (see Definition \ref{strong}).
By Lemma \ref{dec} there is a weak semiorthogonal decomposition
$\db{\gr-A_{\ge 0}}=\langle \cP_{\ge 0}, \T_0 \rangle,$ where
$\T_0$ is equivalent to $\dsinggr{A}.$ Moreover, by formula
(\ref{tidec}) we have the following semiorthogonal decomposition
for $\T_0:$
$$
\T_0=\langle  \kk(0),\dots,  \kk(a+1)\rangle.
$$
Denote by  $E_i,$ where $i=0, \dots, -a-1,$ the modules
$A(i+a+1)/A(i+a+1)_{\ge a}.$ These modules belong to $\T_0$ and
form a full exceptional collection
$$
\T_0=\langle  E_0,\dots,  E_{-(a+1)}\rangle.
$$
Furthermore, this collection is strong, and the algebra of
homomorphisms of this collection  coincides with the algebra
$\qQ(A)=\End_{\gr-A}\left(\bigoplus\limits_{i=a+1}^{0} A(i)\right).$ As
in the previous proposition, consider the object
$E=\bigoplus\limits_{i=0}^{-(a+1)} E_i$ and  the functor
$$
\bR\Hom(E, -):\T_0=\dsinggr{A}\lto \db{\mod-\qQ(A)}.
$$
Again the standard reasoning from \cite{Bon,BK} shows that the
functor $\bR\Hom(E, -)$ is an equivalence of triangulated
categories.
\end{dok}

\begin{example} {\rm
The simplest example here is $A=\kk[x]/x^{n+1}.$ In this case the
triangulated category of singularities $\dsinggr{A}$ has a full
exceptional collection and is equivalent to the bounded derived
category of finite-dimensional representations of the Dynkin quiver
of type $A_n:
\underset{n}{\underbrace{\bullet-\bullet-\cdots-\bullet}},$ because
in this case the algebra $\qQ(A)$ is isomorphic to the path algebra
of this Dynkin quiver. This example is considered in detail in the
paper \cite{Ta}. }
\end{example}

\begin{remark}{\rm
There are other cases in which the triangulated category of
singularities $\dsinggr{A}$ has a full exceptional collection. It
follows from Theorem \ref{mai} that if $a\le 0$ and the derived
category $\db{\qgr A}$ has a full exceptional collection, then
$\dsinggr{A}$ has a full exceptional collection as well. It
happens, for example, in the case that the algebra $A$ is related
to a weighted projective line, an orbifold over $\PP^1$
(see, e.g., \cite{GL}).
}
\end{remark}

\subsection{Categories of Coherent Sheaves for Gorenstein Schemes.}

Let $X$ be a connected projective Gorenstein scheme of dimension $n$
and let $\L$ be a very ample line bundle such that the dualizing
sheaf $\omega_X$ is isomorphic to $\L^{-r}$ for some $r\in\ZZ.$
Denote by $A$ the graded coordinate algebra $\bigoplus_{i\ge
0}H^0(X, \L^{i}).$ The famous Serre theorem \cite{Ser} asserts that
the abelian category of coherent sheaves $\coh(X)$ is equivalent to
the quotient category $\qgr A.$

Assume also that $H^j(X, \L^k)=0$ for all $k\in\ZZ$ when $j\ne 0, n.$
For example, if $X$ is a complete intersection in $\PP^N$ then it
satisfies these conditions.  In this case, Theorem \ref{mai} allows us
to compare the triangulated category of singularities $\dsinggr{A}$
with the bounded derived category of coherent sheaves $\db{\coh(X)}.$
To apply that theorem we need the following lemma.

\begin{lemma}\label{gorcoor}
Let $X$ be a connected projective Gorenstein scheme of dimension
$n.$ Let $\L$ be a very ample line bundle such that $\omega_X\cong
\L^{-r}$ for some $r\in\ZZ$ and $H^j(X, \L^k)=0$ for all $k\in\ZZ$
when $j\ne 0, n.$ Then the algebra $A=\mathop\bigoplus_{i\ge
0}H^0(X, \L^{i})$ is Gorenstein with Gorenstein parameter $a=r.$
\end{lemma}
\begin{dok}
Consider the projection functor $\varPi:\Gr-A\to \QGr A$ and its right
adjoint
$\varOmega: \QGr A\to \Gr-A,$ which is given by the formula
(\ref{funad})
$$
\varOmega\varPi N\cong \bigoplus_{n=-\infty}^{\infty}
\Hom_{\QGr}(\varPi A, \varPi N(n)).
$$
The functor $\varOmega$ has a right derived $\bR\varOmega$ that
is given by the formula
$$
\bR^j\varOmega(\varPi N)\cong \bigoplus_{n=-\infty}^{\infty}
\Ext^{j}_{\QGr}(\varPi A, \varPi N(n))
$$
(see, e.g., \cite[Prop.~7.2]{AZ}).  The assumptions on $X$ and $\L$
imply that $\bR^j\varOmega(\varPi A)\cong 0$ for all $j\ne 0, n.$
Moreover, since $X$ is Gorenstein and $\omega_X\cong \L^{-r},$ Serre
duality for $X$ yields that
$$
\bR^0\varOmega(\varPi A)\cong \bigoplus_{i=-\infty}^{\infty}H^0(X,
\L^i)\cong A
\qquad
\text{ and}
\qquad
\bR^n\varOmega(\varPi A)\cong \bigoplus_{i=-\infty}^{\infty}H^n(X,
\L^i)\cong A^*(r),
$$
where $A^*=\Hom_{\kk}(A, \kk).$
Since $X$ is irreducible, the algebra $A$ is connected.
Since $\varPi$ and $\bR\varOmega$ are adjoint functors, we have
$$
\bR\Hom_{\Gr}(\kk(s), \bR\varOmega(\varPi A))\cong
\bR\Hom_{\QGr}(\varPi\kk(s), \varPi A)=0
$$
for all $s.$ Furthermore, we know that $\bR\Hom_{A}(\kk,
A^{*})\cong\bR\Hom_A(A, \kk)\cong \kk.$ This implies that
$\bR\Hom_A(\kk, A)\cong \kk(r)[-n-1].$ This isomorphism implies that
the affine cone $\Spec A$ is Gorenstein at the vertex, and the
assumption on $X$ now implies that $\Spec A$ is a Gorenstein scheme
(\cite[V, \S 9,10]{Ha}).  Since $\Spec A$ has a finite Krull
dimension, the algebra $A$ is a dualizing complex for itself, i.e.,
it has a finite injective dimension.  Thus, the algebra $A$ is
Gorenstein with parameter $r.$
\end{dok}

\begin{theorem}\label{mai2}
Let $X$ be a connected projective Gorenstein scheme of dimension
$n.$ Let $\L$ be a very ample line bundle such that
$\omega_X\cong\L^{-r}$ for some $r\in\ZZ.$ Suppose $H^j(X, \L^k)=0$
for~all $k\in\ZZ$ when $j\ne 0, n.$ Set $A
:=\mathop\bigoplus_{i\ge 0} H^0(X, \L^i).$ Then the derived
category of coherent sheaves $\db{\coh(X)}$ and the triangulated
category of singularities $\dsinggr{A}$ are related as follows:
\begin{enumerate}
\item[{\rm (i)}]
if $r>0,$ i.e., if $X$ is a Fano variety, then there is a
semiorthogonal decomposition
$$
\db{\coh(X)}=\langle \L^{-r+1},\dots, \O_X, \dsinggr{A}\rangle;
$$
\item[{\rm (ii)}]
if $r<0,$ i.e., if $X$ is a variety of general type, then there is a
semiorthogonal decomposition
$$
\dsinggr{A}=\langle q\kk(r+1),\dots, q\kk,  \db{\coh(X)}\rangle,
$$
where $q:\db{\gr-A}\lto\dsinggr{A}$ is the natural projection;
\item[{\rm (iii)}]
if $r=0,$ i.e., if $X$ is a Calabi-Yau variety, then  there is an
equivalence
$$\dsinggr{A}\stackrel{\sim}{\lto} \db{\coh(X)}.$$
\end{enumerate}
\end{theorem}
\begin{dok}
Since $\L$ is very ample, Serre's theorem implies that the bounded
derived category $\db{\coh(X)}$ is equivalent to the category
$\db{\qgr A},$ where $A=\bigoplus_{i\ge 0}H^0(X, \L^i).$  Since
$H^j(X, \L^k)=0$ for $j\ne 0, n$ and all $k\in\ZZ,$ Lemma
\ref{gorcoor} implies that $A$ is Gorenstein.  Now the theorem
immediately follows from Theorem \ref{mai}.
\end{dok}

\begin{corollary}
Let $X$ be an irreducible projective Gorenstein Fano variety of
dimension $n$ with at most rational singularities. Let $\L$ be a
very ample line bundle such that $\omega_X^{-1}\cong \L^{r}$ for
some $r\in \NN.$ Set $A=\bigoplus_{i\ge 0} H^0(X, \L^i).$ Then the
category $\db{\coh(X)}$ admits a semiorthogonal decomposition of the
form
$$
\db{\coh(X)}=\langle \L^{-r+1},\dots, \O_X, \dsinggr{A}\rangle.
$$
\end{corollary}
\begin{dok}
The Kawamata-Viehweg vanishing theorem (see,
e.g., \cite[Th.~1.2.5]{KMM}) yields  $H^j(X, \L^{k})=0$ for
$j\ne 0, n$ and all $k.$ Hence, we can apply Theorem \ref{mai2}(i).
\end{dok}

\begin{corollary}
Let $X$ be a Calabi-Yau variety. That is, $X$ is an irreducible
projective variety with at most rational singularities, with trivial
canonical sheaf $\omega_X\cong\O_X$ and such that $H^j(X, \O_X)=0$
for $j\ne 0, n.$  Let $\L$ be  some very ample line bundle on $X.$
Set $A=\bigoplus_{i\ge 0} H^0(X, \L^i).$ Then there is an
equivalence
$$
\db{\coh(X)}\cong\dsinggr{A}.
$$
\end{corollary}
\begin{dok}
The variety $X$ has rational singularities hence it is
Cohen-Macaulay. Moreover, $X$ is Gorenstein, because
$\omega_X\cong\O_X.$ The Kawamata-Viehweg vanishing theorem
(\cite[Th.~1.2.5]{KMM}) yields $H^j(X, \L^k)=0$ for $j\ne 0, n$ and
all $k\ne 0.$ Since by assumption $H^j(X, \O_X)=0$ for $j\ne 0, n,$
we can apply Theorem \ref{mai2} (iii).
\end{dok}

\begin{proposition}
Let $X\subset \PP^N$ be a complete intersection of $m$ hypersurfaces
$D_1, \dots, D_m$ of degrees $d_{1}, \dots, d_{m}$ respectively.
Then $X$ and $\L=\O_X(1)$ satisfy the conditions of Theorem
\ref{mai2} with Gorenstein parameter $r=N+1-\sum_{i=1}^{m} d_i.$
\end{proposition}
\begin{dok}
Since the variety $X$ is a complete intersection, it is Gorenstein.
The canonical class $\omega_X$ is isomorphic to $\O(\sum d_i -N-1).$
It can be easily proved by induction on $m$ that $H^j(X, \O_X(k))=0$
for all $k$ and $j\ne 0,n,$ where $n=N-m$ is the dimension of $X.$
Indeed, The base of the induction is clear. For the induction step,
assume that for $Y=D_1\cap\cdots\cap D_{m-1}$ these conditions hold.
Then, consider the short exact sequence
$$
0\lto \O_Y(k-d_m)\lto \O_Y(k)\lto \O_X(k)\lto 0.
$$
Since the cohomologies $H^j(Y, \O_Y(k))$ are $0$ for all $k$ and $j\ne 0,
n+1,$ we obtain that $H^j(X, \O_X(k))=0$ for all $k$ and $j\ne 0, n.$
\end{dok}

Theorem \ref{mai2} can be extended to the case of quotient stacks.
To do this we will need an appropriate generalization of Serre's
theorem \cite{Ser}.  The usual Serre theorem says that if a
commutative connected graded algebra $A=\mathop{\bigoplus}_{i\ge 0}
A_{i}$ is generated by its first component, then the category $\qgr
A$ is equivalent to the category of coherent sheaves $\coh(X)$ on
the projective variety $X=\Proj A.$ (Such equivalence holds for the
categories of quasicoherent sheaves $\Qcoh(X)$ and $\QGr A$ too.)

Consider now  a commutative connected graded $\kk$\!-algebra
$A=\mathop\bigoplus_{i\ge 0} A_i$ that is not necessary generated
by its first component.
The grading on $A$ induces an action of the group $\kk^*$ on the affine
scheme $\Spec A.$ Let $\bz$ be the closed point of $\Spec A$ that
corresponds to the ideal $A_+=A_{\ge 1}\subset A.$ This point is invariant
under the action.

Denote by $\PROJ A$ the quotient stack $\left[(\Spec A\backslash
\bz)\big/\kk^*\right].$
(Note that there is a natural map $\PROJ A \to \Proj A,$ which is an
isomorphism if the algebra $A$ is generated by $A_1.$)
\begin{proposition}{\rm (see also \cite{AKO})}\label{serst}
Let $A=\bigoplus_{i\ge 0} A_{i}$ be a
 connected graded finitely generated algebra. Then the category of
(quasi)coherent sheaves on the quotient stack $\PROJ(A)$ is
equivalent to the  category $\qgr A$ (respectively $\QGr A$).
\end{proposition}
\begin{dok}
Let $\bz$ be the closed point on the affine scheme $\Spec A$ that
corresponds to the maximal ideal $A_{+}\subset A.$ Denote by $U$
the complement $\Spec A\backslash\bz.$ We know that the category of
(quasi)coherent sheaves on the stack $\PROJ A$ is equivalent to
the category of $\kk^*$-equivariant (quasi)coherent sheaves on
$U.$ The category of (quasi)coherent sheaves on $U$ is equivalent
to the quotient of the category of (quasi)coherent sheaves on
$\Spec A$ by the subcategory of (quasi)coherent sheaves with
support on $\bz$ (see \cite{Gab}). This is also true for the categories of
$\kk^*$-equivariant sheaves. But the category of (quasi)coherent
$\kk^*$-equivariant sheaves on $\Spec A$ is just the
category $\gr-A$ (resp. $\Gr-A$) of graded modules over $A,$ and the
subcategory of (quasi)coherent sheaves with support on $\bz$
coincides with the subcategory $\tors-A$ (resp.\ $\Tors-A$). Thus, we
obtain that $\coh(\PROJ A)$ is equivalent to the quotient category
$\qgr A=\gr-A/\tors-A$ (and $\Qcoh(\PROJ A)$ is equivalent to
$\QGr A=\Gr-A/\Tors-A$).
\end{dok}

\begin{corollary}
Assume that the Noetherian Gorenstein connected graded algebra $A$
from Theorem \ref{mai} is finitely generated and commutative. Then
in place of the bounded derived category $\db{\qgr A}$ in Theorem
\ref{mai} we can substitute the category $\db{\coh(\PROJ A)},$ where
$\PROJ A$ the quotient stack $\left[(\Spec A\backslash
\bz)\big/\kk^*\right].$
\end{corollary}

\section{Categories of Graded D-branes of Type B in Landau-Ginzburg Models}

\subsection{Categories of Graded Pairs}\label{brany}

Let $B=\bigoplus_{i\ge 0} B_{i}$ be a finitely generated
connected graded algebra over a field $\kk.$
Let $W\in B_n$ be a central element of degree $n$
that is not a zero-divisor,
i.e., $Wb=bW$ for any $b\in B$ and $bW=0$ only for $b=0.$
Denote by $J$ the two-sided ideal $J := WB=BW$
and  denote by $A$ the quotient graded algebra $B/J.$

With any such element $W\in B_n$ we can associate two categories:
an exact category $\gabe(W)$ and a triangulated category
$\gdt(W).\footnote{One can also construct a  differential graded category
whose homotopy category  is equivalent to $\gdt.$}$
Objects of  these categories are  ordered pairs
$$
\ove{P}:=\Bigl( \xymatrix{ P_1 \ar@<0.6ex>[r]^{p_1} &P_0
\ar@<0.6ex>[l]^{p_0} } \Bigl),
$$
where $P_0, P_1\in \gr-B$ are finitely generated free graded right
$B$\!-modules, $p_1$ is a map of degree $0,$ and $p_0$ is a map of
degree $n$ (i.e a map from $P_0$ to $P_1(n)$) such that the
compositions ${p_0 p_1}$ and ${p_1(n) p_0}$ are the left
multiplications by the element $W.$ A morphism $f:\ove{P}\to\ove{Q}$
in the category $\gabe(W)$ is a pair of morphisms $f_1: P_1\to Q_1$
and $f_0: P_0\to Q_0$ of degree $0$ such that $f_1(n)p_0=q_0 f_0$
and $q_1 f_1=f_0 p_1.$ The morphism $f=(f_1, f_0)$ is null-homotopic
if there are two morphisms $s: P_0\to Q_1$ and $t:P_1\to Q_0(-n)$
such that $f_1=q_0(n) t + s p_1$ and $f_0=t(n) p_0 + q_1 s.$
Morphisms in the category $\gdt(W)$ are the classes of morphisms in
$\gabe(W)$ modulo null-homotopic morphisms.

In other words, objects of both categories are quasi-periodic infinite
sequences
$$
\ul{K}^{\Ddot}:=\{\cdots\lto K^{i}\stackrel{k^i}{\lto} K^{i+1}
\stackrel{k^{i+1}}{\lto} K^{i+2}\lto\cdots\}
$$
of morphisms in $\gr-B$ of  {\it free} graded right $B$\!-modules such
that the composition $k^{i+1}k^i$  of any two consecutive morphisms
is equal to multiplication by $W.$ The quasi-periodicity  property
here means that $\ul{K}^{\Ddot}[2]= \ul{K}^{\Ddot}(n).$ In
particular,
$$
K^{2i-1}\cong P_1(i\cdot n),\;K^{2i}\cong P_{0}(i\cdot n),\;
k^{2i-1}=p_1(i\cdot n),\; k^{2i}=p_0(i\cdot n).
$$
A morphism $f: \ul{K}^{\Ddot}\lto \ul{L}^{\Ddot}$ in the category
$\gabe(W)$
is a family of morphisms $f^{i}: K^{i}\lto L^{i}$ in $\gr-B$
that is quasi-periodic, i.e $f^{i+2}=f^{i}(n),$
and that commutes with $k^{i}$ and $l^i,$ i.e., $f^{i+1} k^{i}=l^{i} f^{i}.$

Morphisms in the category $\gdt(W)$ are morphisms in $\gabe(W)$ modulo
null-homotopic morphisms, and we consider only quasi-periodic
homotopies, i.e.,  families $s^i: K^i\lto L^{i-1}$ such that
$s^{i+2}=s^i(n).$

\begin{definition}\label{grdbr}
The category $\gdt(W)$ constructed above will be called {\sf the
  category of graded D-branes of type
B} for the pair $(B=\bigoplus_{i\ge 0} B_i, W).$
\end{definition}
\begin{remark}
{\rm If $B$ is commutative, then we can consider the affine scheme
$\Spec B.$ The grading on $B$ corresponds to an action of the algebraic group
$\kk^*$ on $\Spec B.$ The element $W$ can be viewed as a regular
function on $\Spec B$ that is semi-invariant with respect to this
action. This way, we get a singular Landau-Ginzburg model $(\Spec
B, W)$ with an action of the torus $\kk^*.$ Thus, Definition \ref{grdbr}
is a definition of the category of {\it graded} D-branes of type B for
this model (see also \cite{HW,Wa}).  }
\end{remark}

It is clear that the category $\gabe(W)$ is an exact category (see
\cite{Qu} for the definition) with monomorphisms and epimorphisms
being the componentwise monomorphisms and epimorphisms. The
category $\gdt(W)$ can be endowed with a natural structure of a
triangulated category. To exhibit this structure we have to define a
translation functor $[1]$ and a class of exact triangles.

The translation functor is usually defined as a functor that takes an
object $\ul{K}^{\Ddot}$ to the object
$\ul{K}^{\Ddot}[1],$ where $K[1]^i=K^{i+1}$ and $d[1]^i=-d^{i+1},$
and takes a morphism $f$ to the morphism
$f[1],$ which coincides  with $f$ componentwise.

For any
morphism $f: \ul{K}^{\Ddot}\to \ul{L}^{\Ddot}$ from the category
$\gabe(W)$
we define a mapping cone $\ul{C}^{\Ddot}(f)$ as an object
$$
\ul{C}^{\Ddot}(f)=\{\cdots \lto L^i\oplus K^{i+1} \stackrel{c^i}{\lto}
L^{i+1}\oplus K^{i+2} \stackrel{c^{i+1}}{\lto}L^{i+2}\oplus
K^{i+3}\lto\cdots\}
$$
such that
$$
\qquad c^i=
\begin{pmatrix}
l^i & \hphantom{-}f^{i+1}\\
0 & -k^{i+1}
\end{pmatrix}.
$$
There are  maps $g: \ul{L}^{\Ddot}\to \ul{C}^{\Ddot}(f), \; g=(\id,
0)$ and $h: \ul{C}^{\Ddot}(f)\to \ul{K}^{\Ddot}[1],\; h=(0, -\id).$

Now we define a standard triangle in the category $\gdt(W)$
as a triangle of the form
$$
\ul{K}^{\Ddot}\stackrel{f}{\lto} \ul{L}^{\Ddot}\stackrel{g}{\lto}
\ul{C}^{\Ddot}(f)
\stackrel{h}{\lto} \ul{K}^{\Ddot}[1]
$$
for some $f\in \gabe(W).$

\begin{definition}
A triangle $\ul{K}^{\Ddot}{\to}\ul{L}^{\Ddot}{\to}\ul{M}^{\Ddot}
{\to}\ul{K}^{\Ddot}[1]$ in
$\gdt(W)$ will be  called  an exact (distinguished)  triangle if it is
isomorphic to a standard triangle.
\end{definition}

\begin{proposition}\label{trstr}
The category $\gdt(W)$ endowed with the translation functor
$[1]$ and the above class of exact triangles becomes a
triangulated category.
\end{proposition}
We omit the proof of this proposition, which is more or less the same as
the proof of the analogous result for a usual homotopic category (see,
e.g., \cite{GM}).

\subsection{Categories of Graded Pairs and Categories of Singularities}

With any object $\ul{K}^{\Ddot}$ as above, one associates a short
exact
sequence
\begin{equation}\label{shseq}
0\lto K^{-1} \stackrel{k^{-1}}{\lto} K^0 \lto \Coker k^{-1} \lto 0.
\end{equation}

We can attach to an object $\ul{K}^{\Ddot}$ the right $B$\!-module
$\Coker k^{-1}.$  It can be easily checked that the multiplication
by $W$ annihilates it.  Hence, the module $\Coker k^{-1}$ is
naturally a right $A$\!-module, where $A=B/J$ with $J = WB=BW.$ Any
morphism $f: \ul{K}^{\Ddot}\to \ul{L}^{\Ddot}$ in $\gabe(W)$ induces
a morphism between cokernels. This construction defines a functor
$\Cok: \gabe(W)\lto \gr-A.$ Using the functor $\Cok$ we can
construct an exact functor between triangulated categories $\gdt(W)$
and $\dsinggr{A}.$
\begin{proposition}\label{funcF} There is a functor $F$ that completes
the following commutative diagram:
$$
\begin{CD}
\gabe(W) @>\Cok>> & \gr-A\\
@VVV&@VVV\\
\gdt(W) @>F>>& \dsinggr{A}.
\end{CD}
$$
Moreover, the functor $F$ is an exact functor between triangulated
categories.
\end{proposition}
\begin{dok}
We have the functor  $\gabe(W)\lto \dsinggr{A},$ which is the
composition of $\Cok$ and the natural functor from $\gr-A$ to
$\dsinggr{A}.$ To prove the existence of a functor $F$ we need to
show that any morphism ${f}: \ul{K}^{\Ddot} \to \ul{L}^{\Ddot}$
that is null-homotopic goes to the $0$\!-morphism in $\dsinggr{A}.$
Fix a homotopy $s=(s^i)$ with $s^i: K^i\to L^{i-1}.$ Consider the
following decomposition of ${f}:$
$$
\xymatrix{ K^{-1} \ar[d]_{(s^{-1}, f^{-1})} \ar[r]^{k^{-1}}  & K^0
\ar[d]^{(s^0, f^0)}\ar[r] & \Coker k^{-1} \ar[d]
\\
L^{-2}\oplus L^{-1} \ar[d]_{pr} \ar[r]^{u^{-1}} & L^{-1}\oplus
L^0\ar[d]^{pr}
\ar[r] & L^0\otimes_{B} A \ar[d]
 & {}\save[]*\txt{where}\restore&
{u^{-1}=
\begin{pmatrix}
-l^{-2} & \id\\
0 & l^{-1}
\end{pmatrix},
}
\\
L^{-1} \ar[r]^{l^{-1}} & L^0
\ar[r]
& \Coker
l^{-1} }
$$ This yields a decomposition of $F({f})$ through a locally free
object $L^0\otimes_{B} A.$ Hence, $F({f})=0$ in the category
$\dsinggr{A}.$  By Lemma \ref{extv}, which is proved below, the
tensor product $\ul{K}^{\Ddot}\otimes_{B}A$ is an acyclic complex.
Hence, there is an exact sequence $ 0\to\Coker k^{-1}\to
K^{1}\otimes_B A\to \Coker k^{0}\to 0.  $ Since $K^{1}\otimes_B A$
is free, we have $\Coker k^{0}\cong \Coker k^{-1}[1]$ in
$\dsinggr{A}.$ But $\Coker k^{0}=F(\ul{K}^{\Ddot}[1]).$ Hence, the
functor $F$ commutes with translation functors. It is easy to see
that $F$ takes a standard triangle in $\gdt(W)$ to an exact triangle
in $\dsinggr{A}.$ Thus, $F$ is exact.
\end{dok}

\begin{lemma}
The functor $\Cok$ is full.
\end{lemma}
\begin{dok}
Any map $g: \Coker k^{-1} \to \Coker l^{-1}$ between $A$\!-modules can
be considered as the map of $B$\!-modules and can be extended to a map
of short exact sequences
$$
\begin{CD}
0 @>>> K^{-1} @>k^{-1}>> K^0 @>>>  \Coker k^{-1} @>>> 0\\
&& @V{f^{-1}}VV @VV{f^0}V @VV{g}V \\
0 @>>> L^{-1} @>l^{-1}>> L^0 @>>>  \Coker l^{-1} @>>> 0,
\end{CD}
$$
because $K^0$ is free. This gives us a sequence
of morphisms $f=(f^{i}), i\in\ZZ,$ where $f^{2i}=f^0(in)$ and
$f^{2i-1}=f^{-1}(in).$
To prove the lemma it is
sufficient  to show that the family $f$ is a map from
$\ul{K}^{\Ddot}$ to $\ul{L}^{\Ddot},$ i.e
$f^1 k^0=l^0 f^0.$ Consider the sequence of equalities
$$
l^1(f^1 k^0 - l^0f^0)= f^2k^1 k^0 - W f^0 = f^2 W - W f^0 = f^0(2) W -
W f^0= 0.
$$
Since $l^1$ is an embedding, we obtain that $f^1 k^0= l^0 f^0.$
\end{dok}

\begin{lemma}\label{extv}
For any sequence $\ul{K}^{\Ddot}\in \gabe(W)$
the tensor product $\ul{K}^{\Ddot}\otimes_{B}A$ is an acyclic complex
of $A$\!-modules
and the $A$\!-module $\Coker k^{-1}$
satisfies the condition
$$
\Ext^{i}_A(\Coker k^{-1}, A)=0
\qquad
\text{for all}
\quad
i>0.
$$
\end{lemma}
\begin{dok}
It is clear that $\ul{K}^{\Ddot}\otimes_B A$ is a complex.
Applying the snake lemma to the commutative diagram
$$
\begin{CD}
0@>>> K^{i-2} @>k^{i-2}>> K^{i-1} @>>> \Coker k^{i-2}@>>> 0\\
&&@V{W}VV  @VV{W}V  @VV{0}V\\
0@>>> K^{i} @>k^i>>  K^{i+1} @>>>  \Coker k^{i}@>>> 0,\\
\end{CD}
$$
we obtain an exact sequence
$$
0\to \Coker k^{i-2}\lto K^i\otimes_{B} A \stackrel{k_i|_{W}}{\lto} K^{i+1}\otimes_B A \lto
\Coker k^i \to 0.
$$
This implies that $\ul{K}^{\Ddot}\otimes_B A$ is an acyclic complex.

Further, consider the dual sequence of left $B$\!-modules
$\ul{K}^{\Ddot\vee},$ where $\ul{K}^{\Ddot\vee}\cong
\Hom_B(\ul{K}^{\Ddot}, B).$ For the same reasons as above,
$A\otimes_{B} \ul{K}^{\Ddot\vee}$ is an acyclic complex.  On the
other hand, the cohomologies of the complex $ \{(K^{0})^{\vee}\lto
(K^{-1})^{\vee}\lto (K^{-2})^{\vee}\lto\cdots\} $ are isomorphic to
$\Ext^{i}_A(\Coker k^{-1}, A).$ And so, by the acyclicity of
$A\otimes_{B} \ul{K}^{\Ddot\vee},$ they are equal to $0$ for all
$i>0.$
\end{dok}

\begin{lemma}\label{embob}
If $F \ul{K}^{\Ddot}\cong 0,$ then $\ul{K}^{\Ddot}\cong 0$ in
$\gdt(W).$
\end{lemma}
\begin{dok}
If $F \ul{K}^{\Ddot}\cong 0,$ then the $A$\!-module $\Coker k^{-1}$
is perfect as a complex
of $A$\!-modules.
Let us show that $\Coker k^{-1}$ is projective in this case.
Indeed, there is a natural number $m$ such that $\Ext^{i}_A(\Coker
k^{-1}, N)=0$
for any $A$\!-module $N$ and any $i\ge m.$
Considering the exact sequence
$$
0\to \Coker k^{-2m-1}\to K^{-2m}\otimes_B A\to\cdots\to
K^{-1}\otimes_B A\to K^0\otimes_B A\to \Coker k^{-1}\to 0
$$
and taking into account that all $A$\!-modules $K^i\otimes_B A$ are
free, we find that
for all modules $N$
$\Ext^{i}_A(\Coker k^{-2m-1}, N)=0$ when $i>0.$
Hence, $\Coker k^{-2m-1}$ is a projective $A$\!-module.
This implies that $\Coker k^{-1}$ is also projective, because it is
isomorphic to
$\Coker k^{-2m-1}(-mn).$

Since $\Coker k^{-1}$ is projective, there is a map $f: \Coker k^{-1}
\to K^{0}\otimes_B A$ that splits the epimorphism $\pr:
K^{0}\otimes_{B} A\to\Coker k^{-1}.$ It can be lifted to a map from
the complex $\{ K^{-1}\stackrel{k^{-1}}{\lto} K^{0} \}$ to the
complex $\{K^{-2}\stackrel{W}{\lto} K^{0}\}.$ Denote the lift by
$(s^{-1}, u).$ Consider the following diagram:
$$
\begin{CD}
K^{-1} @>k^{-1}>> & K^0 @>>> & \Coker k^{-1}\\
@V{s^{-1}}VV & @VV{u}V & @VV{f}V\\
K^{-2} @>W>> & K^0 @>>> & K^{0}\otimes_B A\\
@V{k^{-2}}VV & @VV{\id}V & @VV{pr}V\\
K^{-1} @> k^{-1}>> & K^0 @>>> & \Coker k^{-1}.
\end{CD}
$$
Since the composition $\pr f$ is identical, the map $(k^{-2} s^{-1}, u)$
from the pair $\{ K^{-1}\stackrel{k^{-1}}{\lto} K^0 \}$ to itself is
homotopic to the identity map. Hence, there is a map
$s^{0}: K^0\to K^{-1}$ such that
$$
\id_{K^{-1}} - k^{-2} s^{-1}= s^{0} k^{-1}
\qquad
\text{and}
\qquad
k^{-1} s^{0}=\id_{K^0}-u.
$$
Moreover, we have the following equalities:
$$
0= (u k^{-1}- W s^{-1})=(u k^{-1}- s^{-1}(n)W)=(u - s^{-1}(n) k^0) k^{-1}.
$$ This gives us that $u= s^{-1}(n) k^0,$ because there are no maps
from $\Coker k^{-1}$ to $K^0.$ Finally, we get the sequence of
morphisms $s^{i}: K^{i}\lto K^{i-1},$ where $s^{2i-1}=s^{-1}(in),
s^{2i}=s^0(in),$ such that $ k^{i-1}s^{i}+ k^{i}s^{i+1}=\id.$ Thus
the identity morphism of the object $\ul{K}^{\Ddot}$ is
null-homotopic. Hence, the object $\ul{K}^{\Ddot}$ is isomorphic to
the zero object in the category $\gdt_{0}(W).$
\end{dok}

\begin{theorem}\label{main2}
The  exact functor $F:\gdt(W)\lto \dsinggr{A}$ is fully faithful.
\end{theorem}
\begin{dok}
By Lemma \ref{extv} we have $\Ext^{i}_A(\Coker k^{-1}, A)=0$ for $i>0.$
Now, Proposition \ref{stab} gives an isomorphism
$$
\Hom_{\dsinggr{A}} (\Coker k^{-1} , \Coker l^{-1})\cong
\Hom_{\gr-A}  (\Coker k^{-1} , \Coker l^{-1})/\R,
$$
where $\R$ is the subspace of morphisms factoring through
projective modules. Since the functor $\Cok$ is full, we get that the functor
$F$ is also full.

Next we  show that $F$ is faithful. The reasoning is standard. Let
$f: \ul{K}^{\Ddot}\to \ul{L}^{\Ddot}$ be a morphism for which
$F(f)=0.$ Include $f$ in an exact triangle $
\ul{K}^{\Ddot}\stackrel{f}{\lto}
\ul{L}^{\Ddot}\stackrel{g}{\lto}\ul{M}^{\Ddot}. $ Then the identity
map of $F \ul{L}^{\Ddot}$ factors through the map $F
\ul{L}^{\Ddot}\stackrel{Fg}{\lto} F \ul{M}^{\Ddot}.$ Since $F$ is
full, there is a map $h: \ul{L}^{\Ddot}\to \ul{L}^{\Ddot}$ factoring
through $g: \ul{L}^{\Ddot}\to \ul{M}^{\Ddot}$ such that $Fh=\id.$
Hence, the cone $\ul{C}^{\Ddot}(h)$ of the map $h$ goes to zero under
the functor $F.$ By Lemma \ref{embob} the object $\ul{C}^{\Ddot}(h)$
is the zero object as well, so $h$ is an isomorphism. Thus $g:
\ul{L}^{\Ddot}\to \ul{M}^{\Ddot}$ is a split monomorphism and $f=0.$
\end{dok}

\begin{theorem}\label{main3}
Suppose that the algebra $B$ has a finite homological dimension.
Then the functor $F:\gdt(W)\lto \dsinggr{A}$ is an equivalence.
\end{theorem}
\begin{dok}
We know that $F$ is fully faithful. To prove the theorem we need to
show that each object $T\in \dsinggr{A}$ is isomorphic to $F
\ul{K}^{\Ddot}$ for some $\ul{K}^{\Ddot}\in\gdt(W).$

The algebra $B$ has a finite homological dimension and as a
consequence, it has a finite injective dimension. This implies that
$A=B/J$ has a finite injective dimension too.  By Lemma \ref{lfint}
any object $T\in \dsinggr{A}$ is isomorphic to the image of an
$A$\!-module $M$ such that $\Ext^i_A(M, A)=0$ for all $i>0.$ This
means that the object $D(M)=\bR\Hom_A(M, A)$ is a left
$A$\!-module. We can consider a projective resolution
$\ul{Q}^{\Ddot}\to D(M).$ The dual of $\ul{Q}^{\Ddot}$ is a right
projective
$A$\!-resolution
$$
0\lto M\lto\{ P^0\lto P^1\lto\cdots\}.
$$
Consider $M$ as $B$\!-module and chose an
epimorphism $K^0\twoheadrightarrow M$ from the free $B$\!-module $K^0.$
Denote by $k^{-1}: K^{-1}\to K^0$ the kernel of this map.

The short exact sequence $0\to B\stackrel{W}{\to}B\to A\to 0$ implies
 that for a projective $A$\!-module $P$ and any $B$\!-module $N$ we
 have equalities $\Ext^i_{B}(P, N)=0$ when $i>1.$ This also yields
 that $\Ext^i_{B}(M, N)=0$ for $i>1$ and any $B$\!-module $N,$
 because $M$ has a right projective $A$\!-resolution and the algebra
 $B$ has finite homological dimension. Therefore, $\Ext^i_B(K^{-1}, N)=0$
 for $i>0$ and any $B$\!-module $N,$ i.e., $B$\!-module $K^{-1}$ is
 projective.  Since $A$ is connected and finitely generated, any
 graded projective module is free. Hence, $K^{-1}$ is free.

Since the multiplication on $W$ gives the zero map on $M,$ there is
a map $k^0: K^0\to K^{-1}(n)$ such that $k^0 k^{-1}=W$ and
$k^{-1}(n)k^0=W.$ This way, we get a sequence $\ul{K}^{\Ddot}$ with
$$
K^{2i}\cong K^0(i\cdot n),\; K^{2i-1}=K^{-1}(i\cdot n),\;
k^{2i}=k^0(i\cdot n),\;
k^{2i-1}=k^{-1}(i\cdot n),
$$
and this sequence is an object of $\gdt(W)$ for which
$F\ul{K}^{\Ddot}\cong T.$
\end{dok}

\subsection{Graded D-branes of Type B and Coherent Sheaves}\label{grdbrB}

By a Landau-Ginzburg model we mean the following data: a smooth
variety $X$ equipped with a symplectic K\"ahler form $\omega,$
closed real 2-form ${\mathcal B},$ which is called a B-field, and a
regular nonconstant function $W$ on $X.$  The function $W$ is called
the superpotential of the Landau-Ginzburg model.  Since for the
definition of D-branes of type B a symplectic form and B-field are
not needed, we do not fix them.

With any point $\lambda\in \AA^1$ we can associate a triangulated
category $\dt_{\lambda}(W).$ We give a construction of these
categories under the condition that $X=\Spec(B)$ is affine (see
\cite{KL, Tr}). The category of coherent sheaves on an affine scheme
$X=\Spec(B)$ is the same as the category of finitely generated
$B$\!-modules.  The objects of the category $\dt_{\lambda}(W)$ are
ordered pairs $ \ove{P}:=\bigl( \xymatrix{ P_1 \ar@<0.6ex>[r]^{p_1}
&P_0 \ar@<0.6ex>[l]^{p_0} } \bigl), $ where $P_0, P_1$ are finitely
generated projective $B$\!-modules and the compositions ${p_0 p_1}$
and ${p_1 p_0}$ are the multiplications by the element
$(W-\lambda)\in B.$  The morphisms in the category $\dt(W)$ are the
classes of morphisms between pairs modulo null-homotopic morphisms,
where a morphism $f:\ove{P}\to\ove{Q}$ between pairs is a pair of
morphisms $f_1: P_1\to Q_1$ and $f_0: P_0\to Q_0$ such that
$f_1p_0=q_0 f_0$ and $q_1f_1=f_0 p_1.$  The morphism $f$ is
null-homotopic if there are two morphisms $s: P_0\to Q_1$ and
$t:P_1\to Q_0$ such that $f_1=q_0 t + s p_1$ and $f_0=t p_0 + q_1
s.$

We define a category of D-branes of type B (B-branes) on $X=\Spec(B)$
with the superpotential $W$ as the product
$\dt(W)=\mathop{\prod}_{\lambda\in \AA^1} \dt_{\lambda}(W).$

It was proved in the paper \cite[Cor.~3.10]{Tr} that the category
$\dt_{\lambda}(W)$ for smooth affine $X$ is equivalent to the
triangulated category of singularities $\dsing{X_{\lambda}},$ where
$X_{\lambda}$ is the fiber over $\lambda\in\AA^1.$ Therefore, the
category of B-branes $\dt(W)$ is equivalent to the product
$\mathop{\prod}_{{\lambda}\in \AA^1} \dsing{X_{\lambda}}.$ For
no affine $X$ the category $\mathop{\prod}_{{\lambda}\in \AA^1}
\dsing{X_{\lambda}}$ can be considered as a definition of the
category of D-branes of type B.  Note that in the affine case,
$X_{\lambda}$ is $\Spec(A_{\lambda}),$ where
$A_{\lambda}=B/(W-\lambda)B,$ and hence the triangulated categories
of singularities $\dsing{X_{\lambda}}$ is the same as the category
$\dsing{A_\lambda}.$

Assume now that there is an action of the group $\kk^*$ on the
Landau-Ginzburg model $(X, W)$ such that the superpotential $W$ is
semi-invariant of  weight $d.$ Thus, $X=\Spec(B)$ and
$B=\bigoplus_{i}B_i$ is a graded algebra. The superpotential $W$ is
an element of $B_d.$ Let us assume that $B$ is positively graded and
connected. In this case, we can consider the triangulated category
of graded B-branes $\gdt(W),$ which was constructed in  Section
\ref{brany} (see Definition \ref{grdbr}).

Denote by $A$ the quotient graded algebra $B/WB.$ We see that the
affine variety $\Spec(A)$ is the fiber $X_0$ of $W$ over the point
$0.$ Denote by $Y$ the quotient stack $\left[(\Spec(A)\setminus
\bz)/\kk^*\right],$ where $\bz$ is the point on $\Spec(A)$
corresponding to the ideal $A_{+}.$ Theorems \ref{mai}, \ref{main3}
and Proposition \ref{serst} allow us to establish a relation between
the triangulated category of graded B-branes $\gdt(W)$ and the bounded
derived category of coherent sheaves on the stack $Y.$

First, Theorem \ref{main3} gives us the equivalence $F$ between the
triangulated category of graded B-branes $\gdt(W)$ and the
triangulated category of singularities $\dsinggr{A}.$ Second,
Theorem \ref{mai} describes the relationship between the category
$\dsinggr{A}$ and the bounded derived category $\db{\qgr A}.$ Third,
the category $\db{\qgr A}$ is equivalent to the derived category
$\db{\coh(Y)}$ by Proposition \ref{serst}.  In the particular case,
that $X$ is the affine space $\AA^N$ with the standard action of the
group $\kk^*,$ we get the following result.
\begin{theorem}\label{main4}
Let $X$ be the affine space $\AA^N$ and let $W$ be a homogeneous
polynomial of degree $d.$ Let $Y\subset \PP^{N-1}$ be the hypersurface
of degree $d$ that is given by the equation $W=0.$ Then, there is the
following relation between the triangulated category of graded
B-branes $\gdt(W)$ and the derived category of coherent sheaves
$\db{\coh(Y)}:$
\begin{enumerate}
\item[{\rm (i)}]
if $d<N,$ i.e., if $Y$ is a Fano variety, there is a semiorthogonal
decomposition
$$
\db{\coh(Y)}=\langle \O_Y(d-N+1),\dots, \O_Y, \gdt(W)\rangle;
$$
\item[{\rm (ii)}]
if $d>N,$ i.e., if $X$ is a variety of general type, there is a
semiorthogonal decomposition
$$
\gdt(W)=\langle F^{-1}q(\kk(r+1)),\dots, F^{-1}q(\kk),  \db{\coh(Y)}\rangle,
$$
where $q:\db{\gr-A}\lto\dsinggr{A}$ is the natural projection, and
$F: \gdt\stackrel{\sim}{\lto}\dsinggr{A}$ is the equivalence constructed
in Proposition \ref{funcF};
\item[{\rm (iii)}]
if $d=N,$ i.e., if $Y$ is a Calabi-Yau variety,  there is an
equivalence
$$\gdt(W)\stackrel{\sim}{\lto} \db{\coh(Y)}.$$
\end{enumerate}
\end{theorem}
\begin{remark}
{\rm We can also consider a weighted action of the torus $\kk^*$ on
the affine space $\AA^N$ with positive weights $(a_1,\dots, a_N),$
$a_i >0$ for all $i.$ If the superpotential $W$ is quasi-homogeneous
then we have the category of graded B-branes $\gdt(W).$ The
polynomial $W$ defines an orbifold (quotient stack)
$Y\subset\PP^{N-1}(a_1,\dots, a_N).$ The orbifold $Y$ is the
quotient of $\Spec(A)\backslash \bz$ by the action of $\kk^*,$ where
$A=\kk[x_1,\dots,x_N]/W.$ Proposition \ref{serst} gives the
equivalence between $\db{\coh(Y)}$ and $\db{\qgr A}.$ And Theorem
\ref{mai} shows that we get an analogue of Theorem \ref{main4} for
the weighted case as well.  }
\end{remark}


\end{document}